\documentclass{amsproc}

\usepackage{amsmath,amsthm,amssymb}
\usepackage{graphicx}
\usepackage{wrapfig}
\usepackage{enumitem}
\usepackage{esvect}
\usepackage{bm}
\usepackage{lipsum} 
\usepackage{mathrsfs}
\usepackage{float}
\usepackage{commath}
\usepackage{mathtools}
\usepackage{ulem}
\usepackage{tabularx}

\usepackage[backref]{hyperref} 
\hypersetup{
hidelinks
}

\newtheorem{theorem}{Theorem}[section]

\newtheorem{lemma}[theorem]{Lemma}

\newtheorem{corollary}[theorem]{Corollary}

\newtheorem{definition}[theorem]{Definition}
\newtheorem{claim}[theorem]{Claim}
\newtheorem{conjecture}[theorem]{Conjecture}

\theoremstyle{definition}

\newtheorem{remark}[theorem]{Remark}
\newtheorem{exercise}[theorem]{Exercise}

\newcommand{\PP}{\mathbb{P}}

\newcommand{\RR}{\mathbb{R}}
\newcommand{\CC}{\mathbb{C}}

\newcommand{\R}{\mathbb{R}}
\newcommand{\C}{\mathbb{C}}

\newcommand{\defeq}{\vcentcolon=}

\def\M{\calM}

\def\calH{\mathcal{H}}

\def\i{\sqrt{-1}\,}

\long\def\frame#1#2#3#4{\hbox{\vbox{\hrule height#1pt
 \hbox{\vrule width#1pt\kern #2pt
 \vbox{\kern #2pt
 \vbox{\hsize #3\noindent #4}
\kern#2pt}
 \kern#2pt\vrule width #1pt}
 \hrule height0pt depth#1pt}}}

\def\beq{\begin{equation}}
\def\eeq{\end{equation}}
\def\bal{\begin{aligned}}
\def\eal{\end{aligned}}
\def\bpf{\begin{proof}}
\def\epf{\end{proof}}
\def\lb#1{\label{#1}}

\def\K{K\"ahler }

\def\ra{\rightarrow}

\def\q{\quad}
\def\vp{\varphi}
\def\i{\sqrt{-1}}
\def\ddbar{\partial\overline{\partial}}

\def\eps{\epsilon}

\newcommand{\bdefn}{\begin{definition}}
\newcommand{\edefn}{\end{definition}}
\newcommand{\bremark}{\begin{remark}}
\newcommand{\eremark}{\end{remark}}
\newcommand{\bexer}{\begin{exercise}}
\newcommand{\eexer}{\end{exercise}}

\def\bconj{\begin{conjecture}}
\def\econj{\end{conjecture}}
\def\bcor{\begin{corollary}}
\def\ecor{\end{corollary}}
\def\bthm{\begin{theorem}}
\def\ethm{\end{theorem}}
\def\blemma{\begin{lemma}}
\def\elemma{\end{lemma}}
\def\blem{\begin{lemma}}
\def\elem{\end{lemma}}
\def\bclaim{\begin{claim}}
\def\eclaim{\end{claim}}
\def\bexer{\begin{exercise}}
\def\eexer{\end{exercise}}

\def\h{\hbox}

\def\DD{\mathbb{D}}
\def\PSH{\mathrm{PSH}}

\def\ddbar{\partial\overline{\partial}}

\def\w{\wedge}

\def\MAeq{Monge--Amp\`ere equation }
\def\q{\quad}

\def\eps{\varepsilon}
\def\ra{\rightarrow}

\makeatletter
\newtheorem*{rep@theorem}{\rep@title}
\newcommand{\newreptheorem}[2]{%
\newenvironment{rep#1}[1]{%
 \def\rep@title{#2 \ref{##1}}%
 \begin{rep@theorem}}%
 {\end{rep@theorem}}}
\makeatother

\newreptheorem{theorem}{Theorem}
\newreptheorem{lemma}{Lemma}
\newreptheorem{proposition}{Proposition}

\PassOptionsToPackage{hyphens}{url}\usepackage{hyperref} 
\hypersetup{
  colorlinks=true,
  citecolor=blue,
  linkcolor=blue, 
  urlcolor=cyan
}

\def\K{K\"ahler }

\def\M{{\mathcal M}}

\def\Cvx{\mathrm{Cvx}}

\long\def\frame#1#2#3#4{\hbox{\vbox{\hrule height#1pt
 \hbox{\vrule width#1pt\kern #2pt
 \vbox{\kern #2pt
 \vbox{\hsize #3\noindent #4}
\kern#2pt}
 \kern#2pt\vrule width #1pt}
 \hrule height0pt depth#1pt}}}

\def\bi#1{\bibitem{#1}}

\begin{document}

\title{Convex meets complex}
\dedicatory{\it In memory of my teacher, Steve Zelditch.}


\author{Yanir A. Rubinstein}
\address{{University of Maryland}}
\email{yanir@alum.mit.edu}
\thanks{Research supported by NSF grants 
DMS-1906370,2204347, BSF grants 2016173,2020329.
Thanks to B. Berndtsson and M. Cwikel for helpful discussions.}

\subjclass[2020]{Primary 52A40, 32Q15, 46B70, 35J96, 32W20, 44A15, 14M25}
\keywords{Convex Geometry, Complex Geometry, Interpolation
of Banach Spaces, Toric Varieties, Legendre Transform, Monge--Amp\`ere equations, Mahler Conjecture}

\date{2024}

\begin{abstract}

Convex geometry and complex geometry have long had fascinating interactions.
This survey offers a tour of a few.

\end{abstract}

\maketitle

\section{Prologue}

A convex function can be obtained from a special type of plurisubharmonic function: 
if $f:\C^n\ra \R$ is a function of $z=x+\i y\in\C^n$ independent of $y\in\R^n$ then
$f$ restricts to a convex function of $x\in\R^n$. 
Similarly, a convex body $K$, i.e., a compact convex set in $\R^n$, 
is the base of a trivial $\R^n$-fibration (also called the {\it tube domain} over $K$)
$T_K\defeq \R^n+ \i K$,
so convex bodies can be viewed as the projections of a special subclass of pseudoconvex domains in $\C^n$.

Naively then, convex geometry and analysis might seem like a sub-field of complex geometry and analysis.
This is of course misleading, since in the study of several complex variables
many techniques concentrate on holomorphic objects which often are far too rigid for the analysis 
of convex bodies, say. As a matter of fact, in the intricate relationship between
convex and complex, it is oftentimes the complex side that benefits from convex results
and techniques.

The purpose of this survey is to give a tour of some results, techniques, and
highlights of this relationship. These are motivated by the author's own 
taste and limited knowledge and do not reflect any sort of ``canonical" choices.
At best, they might serve to entice the reader to explore more.



\section{An early convex--complex result: Study's theorem}

After my post-doctoral year with Steve at Johns Hopkins (2008--2009), I continued
to Stanford and in my first quarter there found myself teaching
the graduate Complex Analysis course (Math 215A). For many years
the course adhered to Ahlfors' text \cite{ahlfors}, but I preferred to experiment
with Gamelin's new (at the time) book that my father was fond of.
I complemented it with several other sources, including, e.g., 
a deep-dive into the technicalities surrounding the Bouligand lemma that was missing from Gamelin's treatment of the uniformization theorem.
For the final 
(take-home) exam, I came across the following beautiful theorem of Study
\cite[p. 109]{Study1913}
that I found in the book of Saks--Zygmund (originally published in 1938
in Polish \cite[p. 224]{SaksZ}).

\begin{theorem}
Let $f$ be a conformal map on $\DD$ whose image is a convex set $C\subset \CC$. 
Let $r\in(0,1)$.
The image of $r\DD=\{z\in\DD\,:\, |z|<r\}$ under $f$~is~a~convex~set.
\end{theorem}

(A conformal map is what we call in higher dimensions a biholomorphic map:
holomorphic and its inverse exists and is holomorphic as well.
Alternatively, this is also classically described as a univalent analytic function
of one complex variable, i.e., an analytic and one-to-one function.)

I do not know if this is the first instance of complex/convex interaction, but
it certainly is an early one, dating back to Study's 1913 book. In fact, Study
discussed related topics in his ICM article a year earlier \cite{Study1912}, 
giving an indication to how innovative and important this must have been
considered at the time. 

It is interesting to note how mathematical statements have evolved over
the past century. Study's original statement reads\footnote{
Wenn in einer gew\"ohnlichen Potenzreihe, die sich nicht auf eine
ganze lineare Funktion ihres Argumentes reduziert, die erste Potenz der
unabh\"angigen Ver\"anderlichen wirklich vorkommt, so geh\"ort zu dieser
Rheihe ein Kreis von endlichem Badius, der mit dem Konvergenzkreis
konzentrisch ist (und unter Umst\"anden mit ihm zusammenfallen wird).

Der Fl\"ache dieses Kreises entspricht ein konvexer Bereich, und
ebenso allen kleineren konzentrischen Kreisfl\"achen, w\"ahrend gr\"oßeren
konzentrischen Kreisfl\"achen konvexe Bereiche nicht mehr entsprechen.
(I am grateful to R. Bamler for the careful translation of the somewhat
archaic German original.)
}
\cite[p. 109]{Study1913}:

\smallskip
\noindent
{\it
If in an ordinary power series, which does not reduce to a whole linear function of its argument, [it is the case that] the first power of the independent variable indeed shows up, then for this series there is a circle of finite radius, which is concentric with the circle of convergence (and which potentially agrees with it). 
The area of this circle corresponds to a convex domain, and in the same way [corresponds to] all smaller concentric areas of circles, while [the] larger concentric circle areas of convex domains correspond to no more [than that].
}

\smallskip

The proof below is, according to Saks--Zygmund, due to Rad\`o.
\bpf
Take $z_0,z_1\in\DD$ with $|z_0|\le|z_1|<r$ and $z_0\ne z_1$. 
Assume first that $z_0=0$. Since $C=f(\DD)\ni\{f(0),f(z)\}$ is convex we know that the line segment
$(1-t)f(0)+tf(z)$ is contained in $C$. Therefore, for each fixed
$t\in[0,1]$, the function
$$
h(z)=f^{-1}\big((1-t)f(0)+tf(z)\big)
$$
is an analytic function of $\DD$ to itself fixing zero.
Schwarz' lemma \cite[p. 260]{Gamelin} implies that $|h(z)|\le|z|$, and putting
$z=z_1$ gives that the line segment
$(1-t)f(0)+tf(z_1)$ is contained in the image of $\{|z|<r\}$, as desired.
For the general case, consider the function
$
h(z)=f^{-1}\big((1-t)f(zz_0/z_1)+tf(z)\big),
$
for which the same argument shows the line segment connecting $f(z_0)$
and $f(z_1)$ is contained in $f(r\DD)$.
\epf

%
%
%
%
%
%
%
%
%

The Schwarz Lemma is proved using the maximum principle. Hence,
the next proof should not come as a surprise (cf. \cite[\S1.3]{ahlfors2}
for the related notion of {\it convex univalent} going back
at least to L\"owner \cite{lowner}).

\bpf[Another proof]
We sidestep boundary regularity technicalities by assuming that 
$f$ extends smoothly to the disc $(1+\eps)\DD$, hence
also to the closed disc $\overline{\DD}$. 
Since $f$ is orientation-preserving, convexity of 
$f(\overline{\DD})=\overline{f(\DD)}=\overline{C}$
means that the slope of the tangent to $\partial C$
is non-decreasing, i.e.,
$$
\bal
\frac{d}{d\theta}\arg\Big[\frac{d}{d\theta} f(re^{\i\theta}) \Big]
&=
\frac{d}{d\theta}\arg\Big[\frac{d}{d\theta} f(re^{\i\theta}) \Big]
\cr
&=
\frac{d}{d\theta}\arg[\i z f'(z)]
\cr
&=
\frac{d}{d\theta}(\pi/2 +\arg z+\arg f'(z))
=
\frac{d}{d\theta}\h{\rm Im}\,(\log z+\log f'(z))
\cr
&=
\h{\rm Re}\,\Big(1+\frac {zf''(z)}{f'(z)}\Big)\ge0
\q \h{for $z\in\partial \DD$}.
\eal
$$
Note that since $f$ is conformal it is invertible and $f'(z)\not=0$
(some care is needed to extend this to the boundary in general).
Since the inequality concerns the real part of an analytic function (hence  harmonic)
it holds on all of $\DD$. In particular this holds on $r\partial \DD$
for any $r\in(0,1)$.
\epf

Study's theorem can be interpreted as the convexity 
of sub-level sets of the Green function of a convex domain, and there are numerous
generalizations to both real and complex higher-dimensions and to
more general underlying PDE
\cite{Gabriel,Kawohl}.

\section{Riesz--Thorin and the complex method of interpolation}
\lb{RTSec}

The Riesz--Thorin theorem in its very simplest form reduces to the H\"older inequality.
It answers the basic question: if $f\in L^p\cap L^q$ for some $0<p<q$ is it true
that also $f\in L^r$ for intermediate $r\in(p,q)$? In fact, it gives precise estimates:
\beq\lb{HolderEq}
||f||_{L^r}\le ||f||_{L^p}^{\frac{\frac1r-\frac1q}{\frac1p-\frac1q}}||f||_{L^q}^{\frac{\frac1p-\frac1r}{\frac1p-\frac1q}}.
\eeq
This does not seem related to convexity at first glance, but it is!
Writing
\beq\lb{lambdaEq}
\lambda:=\frac{\frac1p-\frac1r}{\frac1p-\frac1q},
\eeq
reveals an affine relation of the reciprocals,
$$
\frac1r=\frac{1-\lambda}{p}+\frac{\lambda}{q},
$$
while \eqref{HolderEq} is 
$
\log||f||_{L^r}\le (1-\lambda)\log||f||_{L^p}
+\lambda\log||f||_{L^q}
.
$
Thus, \eqref{HolderEq} is
equivalent to the statement that
$$
F(f,\delta):=\log ||f||_{L^{\frac1\delta}}
$$
is convex in $\delta$, i.e.,
$
F(f,(1-t)\delta_1+t\delta_2)\le (1-t)F(f,\delta_1)+tF(f,\delta_2), \q\forall t\in[0,1].
$

This observation ultimately leads to the first ``grown-up'' version of the Riesz--Thorin theorem,
that we turn to next. For $f\in L^p\cap L^q\cap L^s$, the inequality \eqref{HolderEq} can be rewritten
\beq\lb{Holder2Eq}
\frac{||f||_{L^r}}{||f||_{L^s}}\le 
\bigg(\frac{||f||_{L^p}}{||f||_{L^s}}
\bigg)^{\frac{\frac1r-\frac1q}{\frac1p-\frac1q}}
\bigg(\frac{||f||_{L^q}}{||f||_{L^s}}\bigg)^{\frac{\frac1p-\frac1r}{\frac1p-\frac1q}}.
\eeq
In terms of the operator norm,
$$
||T||_{L^a,L^b}:=\sup_{f\in L^a, Tf\in L^b}\frac{||Tf||_{L^b}}{||f||_{L^a}}.
$$
we thus have 
$
\log||T||_{L^s,L^r}\le
(1-\lambda)||T||_{L^s,L^p}+
\lambda||T||_{L^s,L^q}, \q\hbox{for $\lambda$ satisfying \eqref{lambdaEq}},
$
for $T=Id$, the identity operator, i.e., $\log ||T||_{L^s,L^{\frac1\delta}}$
is convex in $\delta$.
The Riesz--Thorin theorem is then the statement that more generally
$\log ||T||_{L^{\frac1\epsilon},L^{\frac1\delta}}$
is {\it jointly} convex in  $(\epsilon,\delta)$.
Interestingly, Riesz first stated his theorem as such joint convexity%
\footnote{
Convexity was gaining popularity in analysis at that time,
as can be inferred from the very first sentence of Riesz' article:
``Il y a vingt ans, Jensen \'a publi\'e dans les {\it Acta mathematica} un 
m\'emoire important sur les fonctions convexes et les in\'egalit\'es qui s'y rattachent. D\`es lors le r\^ole des fonctions convexes en l'Analyse est devenu de plus en plus manifeste."}
but restricted to linear $T$ and to the triangle $\epsilon,\delta,\epsilon+\delta\in[0,1]$
\cite[Th\'eor\`eme I]{Riesz}. Thorin observed, using complex variables, that
the convexity extends to the square $\epsilon,\delta\in[0,1]$ and is not limited
to linear $T$ \cite{Thorin} (cf. \cite{Peetre}). Thorin's paper is the second convex-complex
interaction I am aware of. The reader is referred to Littlewood's and Peetre's beautiful expositions of this theorem \cite[p. 20]{Littlewood},\cite[\S4.2]{PeetreThorin}. Littlewood particularly raves about Thorin's idea:
{\it ``an idea...the most impudent in mathematics, and brilliantly successful."}

Another way to look at the first version is that we are given two `end-points'
$L^p$ and $L^q$ and we want to `connect' them with a curve of Banach spaces
so that along the curve the convexity properties above hold. This is precisely
the basic idea underlying interpolation theory of Banach spaces. As we saw,
the solution is that the interpolating curve will be 
$L^{r(t)}$ satisfying the second-order ODE
\beq\lb{inverseAffineEq}
\Big(\frac1r\Big)''=0, \q r(0)=p, \q r(1)=q.
\eeq

One may seek to generalize two aspects (at least): one may look at other
PDE, and one may consider interpolation not just between two values but, say,
from a family of Banach spaces parametrized by the circle to the interior of the disc.
In that case it is natural to require subharmonicity in the disc variables
and replace the ODE \eqref{inverseAffineEq} by the Dirichlet
problem for the Laplace equation,
\beq\lb{inverseLaplaceEq}
\Delta\frac1r=0, \q r(\cos\theta,\sin\theta)=p_\theta, \q \theta\in[0,2\pi),
\eeq
where $L^{p_\theta}$ are given Banach spaces parameterized by $S^1$.
This version of the Riesz--Thorin theorem can be stated as follows:
Let $T(\tau):L^{p(\tau)}(X,\mu)\rightarrow L^{q(\tau)}(Y,\nu)$ be a family of operators parametrized by $\tau\in\DD$ and varying analytically in $\tau$.
Suppose $\Delta \frac 1p=0=\Delta \frac 1q$ on $\DD$. Then 
$\log ||T(\tau)||_{L^{p(\tau)}(X,\mu),L^{q(\tau)}(Y,\nu)}$ is {\it subharmonic}.

The key in the proof is to use that there is an analytic function $P(\tau)$
whose real part is $1/p$, as the latter is harmonic (and similarly for $q$)
and then use the mean value characterization of subharmonicity. 

To recover the classical 
version of Riesz--Thorin
one chooses the boundary value $p_\theta$ that is equal to $p_0$
for a set of measure $2\pi(1-t)$ and $p_1$ for the complement of measure $2\pi t$,
and similarly for $q_\theta$. The interpolated $L^{p(0)}$ and $L^{q(0)}$ (at the origin
of $\DD$)
are precisely then ones from \eqref{inverseAffineEq}, since by harmonicity
and the mean value equality
$$
\frac1{p(0)}=\int_0^{2\pi}\frac1{p_\theta}d\theta=\frac{1-t}{p_0}+\frac t{p_1},
$$
and similarly for $q(0)$. Finally, subharmonicity of the operator norms reduces to convexity in $t$ since the mean value inequality reads
$$
\bal
\log ||T(\tau)||_{L^{p(0)}(X,\mu),L^{q(0)}(Y,\nu)}
&\le
\frac1{2\pi}
\int_0^{2\pi}
\log ||T(\tau)||_{L^{p_\theta}(X,\mu),L^{q_\theta}(Y,\nu)}
\cr
&=
\frac1{2\pi}
\Big[
2\pi(1-t)\log ||T(\tau)||_{L^{p_0}(X,\mu),L^{q_0}(Y,\nu)}
\cr
&\q\q\q\q\q
+
2\pi t\log ||T(\tau)||_{L^{p_1}(X,\mu),L^{q_1}(Y,\nu)}
\Big].
\eal$$

This very natural {\it complex} proof of the functional analytic theorem of Riesz--Thorin, in the greater generality of Banach spaces,
goes back to Calder\'on \cite{Calderon0,Calderon} (cf. Lions \cite{lions}, and for an exposition \cite{Sadosky}). In their picture
the interpolating parameter lives on the vertical strip $[0,1]\times \RR$ 
with the endpoint Banach spaces being $(L^{p_0}, L^{q_0})$ on $\{0\}\times \RR$
and $(L^{p_1}, L^{q_1})$ on $\{1\}\times \RR$. Then the interpolating
spaces are also independent of the vertical $\RR$ variable.

Of course, the main point of this theory is that it applies not just to $L^p$ spaces but to rather general Banach spaces. E.g., letting the two endpoints correspond to a convex body and its polar, one obtains as the midpoint the self-dual Euclidean ball of radius 1 (the $L^2$ unit ball), leading to a complex proof of the Santal\'o inequality
\cite{ce,santalo}.

The more general approach described above (over $\DD$)
is due to Cwikel et al. \cite{CCRSW0,CCRSW} as part of their theory of complex interpolation
of {\it families} of Banach spaces. In essence, it removes the $\RR$-independence
and allows bona fide complex 1-dimensional families and not just real 
1-dimensional families as in Calder\'on's setting. Equivalently,
this allows to interpolate families (possibly infinite) of Banach spaces
and not just {\it pairs}.
For example, remarkably, this
turns out to generalize the celebrated Wiener--Masani theorem \cite{WM}, corresponding
to the special case that the spaces on $\partial \DD$ are finite-dimensional
Hilbert spaces, i.e., given by a map $p$ from the circle $\partial \DD$ to the 
symmetric space of symmetric positive-definite matrices. Then
the interpolation family is the unique solution of the Dirichlet problem 
\beq\lb{WMEq}
(P^{-1}P_\tau)_{\bar\tau}=0,\q \tau\in\DD, \q\q P|_{\partial \DD}=p.
\eeq
Reviewing different proofs of the Wiener--Masani theorem would entail (and merits)
a survey in its own right, and we leave that to a future article.

This point of view on interpolation theory has had enormous influence not just
on that theory itself but also on complex geometry. Indeed, 
a decade later Semmes \cite{Semmes1988,S} realized that these 1-parameter
curves (real or complex) can be interpreted as geodesics in an
infinite-dimensional space associated to a (weakly) Riemannian metric. 
The simplest case of Riesz--Thorin simply corresponds to the metric
\beq\lb{oneovermetricEq}
d(L^p,L^q):=|p^{-1}-q^{-1}|
\eeq on $\RR_+$ whose geodesics are \eqref{inverseAffineEq}.
The generalization \eqref{inverseLaplaceEq} no longer corresponds to geodesics,
but rather to harmonic maps from $\DD$ to $\RR_+$ where $\DD$ is equipped
with the standard Euclidean metric and $\RR_+$ with \eqref{oneovermetricEq}.
Of course, geodesics are a special case of harmonic maps \cite{ES}.
Inspired by this 
realization he introduced, in the infinite-dimensional
space of \K metrics \eqref{HLEq}, a notion of geodesics and (twisted) harmonic maps
of the disc. In a nutshell, this involved substituting the Laplace equation \eqref{inverseLaplaceEq} by its higher-dimensional fully nonlinear generalization: the homogeneous complex \MAeq (HCMA). 
The study of Semmes' metric has become a central subfield of \K geometry over the past four decades
(independently discovered by Mabuchi \cite{Mab}, and a decade later by Donaldson \cite{Don})
see, e.g., \cite{BBInv,Ber,Ber-long-geod,BerBerm,Blo,Darvas-AJM,DR,Don2,LT,LV,PS,RWN-IHES,RZ0,RZ1,RZ2,RZ3,SZ1,SZ2,WN},
and the surveys \cite[Chapter 2]{RThesis},\cite{R14,RWN-survey,Darvas-survey,R20}. 
We come back to it in \S\ref{toric-HRMA}.
Coifman--Semmes studied a generalization where the family is parameterized by several complex variables, and turns out to be related to curvature of
vector bundles (such a relation goes back already
to Rochberg \cite{Rochberg} who noticed \eqref{WMEq} is a bundle zero-curvature
equation) and the Yang--Mills equation \cite{CS}
(cf. Slodkowski \cite{slodki,slodki86}). We end this section 
with Semmes' words, that can be a source of inspiration in many fields \cite[p. 156]{Semmes1988}:
{\it ``We must be more adventurous in the theory of interpolation of Banach
spaces."}
Semmes' pioneered the interactions between interpolation theory, several
complex variables and pluripotential theory, PDE, and complex geometry. 
Some ideas and problems that he raised were developed further recently by
Berndtsson--Cordero-Erausquin--Klartag--Rubinstein~\cite{BCKR1,BCKR2},~see~also~\S\ref{CxLegSeg}.

\section{Berezin--Berndtsson--Lagerberg formalism
}
\lb{BBLSec}

\def\dsharp{d^{\#}}

An indispensable tool in complex analysis and geometry is the $\partial\bar\partial$-lemma \cite[p. 149]{GH},
particularly through the integration-by-parts convenience that it 
summons. An illustration of this is furnished by a quick comparative glance
at Aubin's 1984 paper \cite{Aubin1984} that avoided this formalism and later treatments
of the same topic that availed of it (e.g., \cite{Mab2,Tian87}, \cite[\S4.4]{RThesis}).

Berndtsson observed that this same convenience can be, somewhat
artificially, mimicked in the real setting, by adapting an idea of Berezin. 
The gist of it: introduce $n$ dummy (or ``super") variables $\xi\in\RR^n$, and use the 
convenience of positive ``(1,1)-forms" on the ``doubled 
space" agreeing 
formally that $\RR^n$-fibers (parameterized by $\xi$) 
above each point $x\in\RR^n$ have volume 1, i.e., 
$
\int_{\RR}d\xi_i=1, \q i\in\{1,\ldots,n\},
$
or 
\beq
\lb{BerezinEq}
\int_{\RR^n}d\xi_1\w\cdots\w
d\xi_n=1
\eeq
(up to a sign, see \eqref{BerezinModifiedEq} below).
This is what Berezin calls a superintegral \cite[(2.2.1)]{Berezin}.
In Berezin's formalism one has the variables $x$ as here, and then
``anticommuting variables" $\xi_i$. The latter are simply our 
forms $d\xi_i$ that indeed anticommute in the exterior algebra.
Berezin then considers objects he denotes by
$\sum_{\beta} f_\beta(x)\xi^\beta$
(note the dependence only on $x$), 
\cite[(1.3.1)]{Berezin} 
where 
$\beta=(\beta_1,\ldots,\beta_q)$ are multi-indices
and 
$\xi^\beta:=\xi_{\beta_1}\cdots\xi_{\beta_q}$.
These objects in our notation will be $q$-forms 
$\sum_{\beta}f_{\beta}(x) d\xi^\beta$ where
$d\xi^\beta:=d\xi_{\beta_1}\w\cdots\w d\xi_{\beta_p}$.
{\it Berezin's integrals} of these objects are formally integrals over $x$
obtained by (essentially) implicitly multiplying these forms by 
$dx:=dx_1\w\cdots \w dx_n$ and integrating (in $x$ and formally
also in $\xi$, using \eqref{BerezinEq}).
More generally and perhaps more naturally, we will consider objects that we will
think of as ``$(p,q)$-forms''
$\sum_{\alpha,\beta}f_{\alpha,\beta}(x)
dx^\alpha\w d\xi^\beta$ with
$\alpha=(\alpha_1,\ldots,\alpha_p)$ and use the usual
calculus of the exterior algebra of forms on $\RR^{2n}$,
and, when integrating, the Berezin rule \eqref{BerezinEq}.

There are two (more or less equivalent) ways to incorporate derivatives
and the exterior algebra of differential forms and superforms into this formalism.
In the first version \cite{Bern-Arnold}, one  complexifies $\RR^n$ to $\CC^n$
using the standard almost complex structure 
$I_-(\frac{\partial}{\partial x},\frac{\partial}{\partial \xi}):=
(\frac{\partial}{\partial \xi},-\frac{\partial}{\partial x})$
(even though tempting, one should not really think of $I_-$ as acting on the 
coordinates themselves by $(x,\xi)\mapsto(\xi,-x)$).
In the second version \cite{Lagerberg, Larson}, one 
para-complexifies $\RR^n$ to $\RR^n\oplus\RR^n$
instead using the almost 
para-complex structure
$I_+(\frac{\partial}{\partial x},\frac{\partial}{\partial \xi}):=
(\frac{\partial}{\partial \xi},\frac{\partial}{\partial x})$
 \cite{Gadea}.
In both versions, crucially, the bi-degree is determined by
the ``number of $dx$'s and $d\xi$'s", e.g.,
the form
$fdx_1\w d\xi_2$ is a (1,1)-form. Also, in both versions
all objects can only depend on $x$, i.e., $f=f(x)$.
This sort of bi-degree is different from the bi-degree 
associated to the complex  or paracomplex coordinates.
In complex coordinates $dz_j:=dx_j+\i d\xi_j$
with $\overline{dz_j}:=dx_j-\i d\xi_j$
the 
bidegree is determined by counting the number of $dz$'s
and $\overline{dz}$'s. In paracomplex coordinates the same
is true where
$dz_j:=dx_j+\sqrt{1} d\xi_j$ and $\overline{dz_j}:=dx_j-\sqrt{1} d\xi_j$.
Thus, e.g.,
the (real) form $dz_1\w \overline{dz_2}+\overline{dz_1}\w dz_2
=2dx_1\w dx_2-(\sqrt{\pm1})^22d\xi_1\w d\xi_2$ has bidegree (1,1) in the 
usual sense of (para)complex variables, but mixed bidegree (2,0)
and (0,2) in the Berezin--Berndtsson--Lagerberg sense. In essence, we are thinking of the
$x$ coordinates as the ``$dz$"'s and the $\xi$ coordinates as their
conjugates.


The (para)complex structure extends to the exterior algebra of differential forms:
for functions (i.e., 0-forms) $I_\pm f:= f$ and for all other forms
$I_\pm\alpha:=\alpha\circ I_\pm^{-1}$. Note that $I_\pm^{-1}=\pm I_\pm$ on vector
fields. This definition is natural considering the natural pairing
between vector fields and 1-forms. In particular it implies that 
that if $\alpha$ and $X$ are such that $\alpha(X)=1$ then also $(I_\pm\alpha)(I_\pm X)=1$.
In other words, $I_\pm$ commute with the ${}^\sharp$ and ${}^\flat$ operations
of Riemannian geometry (we are using the Euclidean metric). Concretely,
$I_\pm dx_j=d\xi_j$ and $I_\pm d\xi_j=\pm dx_j$.
Thus, if $\alpha$ is a $p$-form, $I_\pm\alpha=(\pm1)^p\alpha\circ I_\pm$, by which 
we mean $(I_\pm\alpha)(Y_1,\ldots,Y_p)=(\pm1)^p\alpha(I_\pm Y_1,\ldots,I_\pm Y_p)$.
The term $(\pm1)^p$ can also be avoided if we instead define, as usual,
$(I_\pm\alpha)(Y_1,\ldots,Y_p)=\alpha(I^{-1}_\pm Y_1,\ldots,I^{-1}_\pm Y_p)$.
Concretely, $I_\pm (dx_1\w dx_2\w d\xi_3)=I_\pm(dx_1)\w I_\pm(dx_2)\w I_\pm(d\xi_3)
=\pm d\xi_1\w d\xi_2\w dx_3$
while $I_\pm (dx_1\w d\xi_2\w d\xi_3)
= d\xi_1\w dx_2\w dx_3$. 
In other words, $I_-$ picks up a minus sign for each appearance of a $d\xi_j$
(while $I_+$ never does). 
The twisted exterior differential $d^{\pm}$ is
$$
d^{\pm}:=I_\pm\circ d\circ I_\pm^{-1},
$$
that when applied to a $p$-form acts by $(\pm 1)^pI_\pm \circ d \circ I_\pm$.
%
%
Because of the previous observations, both of the operators $d^\pm$ actually
have the same formula when acting on forms independent of $\xi$:
%
\begin{align*}
\allowdisplaybreaks
d^\pm
\sum_{\alpha,\beta}f_{\alpha,\beta}(x)
dx^\alpha\w d\xi^\beta
&=\sum_{j=1}^n\frac{\partial}{\partial x^j}\wedge d\xi_j
\Big(
\sum_{\alpha,\beta}f_{\alpha,\beta}(x)
dx^\alpha\w d\xi^\beta
\Big)
\cr
&=
\sum_{\alpha,\beta}
\sum_{j=1}^n\frac{\partial f_{\alpha,\beta}(x)}{\partial x^j}
\wedge d\xi_j\w dx^\alpha\w d\xi^\beta.
\end{align*}

Also, $d^\pm\circ d^\pm=0$ as $I_\pm^2=\pm I$ and $d^2=0$.
When applied to a function $f=f(x)$,
$$
\bal
dd^\pm f
&=
(d\circ I_\pm \circ d\circ I_\pm^{-1})f
=
d\circ I_\pm \circ df
\cr
&=
d\circ I_\pm
\sum_{j=1}^n\frac {\partial f}{\partial x^j}dx^j
=
d
\sum_{j=1}^n\frac {\partial f}{\partial x^j}d\xi^j
=
\sum_{i,j=1}^n\frac {\partial^2 f}{\partial x^i\partial x^j}
dx^i\w d\xi^j,
\eal
$$
that we consider a $(1,1)$-form.
We say this is a positive $(1,1)$-form when the Hessian
of $f$ is a non-negative matrix-valued function. More generally, when
$f$ is convex $dd^{\pm}f$ is a so-called positive $(1,1)$-current.
Lagerberg shows the $dd^{\pm}$-lemma for $d$-closed positive $(1,1)$-currents
\cite[Proposition 2.6]{Lagerberg}.
While that is the most useful case, just as
for the $\partial\bar\partial$-lemma, it holds also in a local
sense for general closed (1,1)-forms \cite[p. 36]{RThesis}.
Note also that $d^-$ is simply the $d^c$ operator familiar to complex analysts
(where, again, we have changed the bi-grading of the differential complex).





A particular example of how this formalism can be useful for writing proofs
that are otherwise a bit tedious is the representation of mixed volumes
in terms of support functions, due to Lagerberg \cite[\S3]{Lagerberg}.
A key observation is that while $\RR^n$ is non-compact and convex bodies
have boundary, there is still
a natural way to do computations much as on compact \K manifolds. 
The inspiration can be seen as coming from toric manifolds and
pluripotential theory. 
In \K geometry, one considers a compact \K manifold $(X,\omega)$.
The \K class $[\omega]\in H^2(M,\RR)$ it determines gives rise to the
(infinite-dimensional)
space $\calH_{\omega}$ of all \K forms representing that class. This notion goes back to
Calabi \cite{Calabi}, who first noticed that the $\ddbar$-lemma
gives that
\beq
\lb{HLEq}
\calH_\omega:=
\{\vp\in C^\infty(X) \,\,:\,\, \omega_\vp:=\omega+\i\ddbar\vp>0\}.
\eeq
This was later generalized to the space of $\omega$-psh function 
$\PSH(X,\omega)\subset L^1(X,\omega^n)$ which corresponds to all \K currents representing $[\omega]$ and subsequently to the space of full-measure $\omega$-psh functions
$\mathcal{E}(X,\omega)\subset \PSH(X,\omega)$ \cite[Definition 1.1]{GZ}.
Lagerberg observed that the latter has a convex analogue, namely the 
infinite-dimensional space
\beq
\lb{HKEq}
\mathcal{E}_K:=
\{f\in \Cvx(\RR^n) \,\,:\,\, \lim_{t\ra\infty}f(t\,\cdot)/t=h_K\}\subset\Cvx(\RR^n),
\eeq
for any $K\in\mathcal{K}_0(\RR^n)$, where $h_K:\RR^n\ra\RR$ is
the support function of $K$, defined by
$h_K(y):=\sup_{x\in K}\langle x,y\rangle$.
The bridge between these two spaces is, as often is the case, 
provided by toric geometry: when $X$ is a toric manifold 
then on the open dense orbit of the $(\CC^*)^n$-action each element
of $\mathcal{E}(X,\omega)$ restricts to an element of
 $\mathcal{E}_P$ where $P$ is an associated Delzant polytope
(unique if we require its barycenter to be the origin), see \S\ref{toricSec}.

On a compact \K manifold 
$\int_X\eta_1\w\cdots\w \eta_n=[\eta_1]\cup\cdots\cup[\eta_n]([X])$ (with the pairing between cohomology (with its cup-product structure)
and homology) that is also denoted as the intersection number 
$[\eta_1].\cdots.[\eta_n]$ (where each $\eta_i\in H^2(M,\RR)$).
In the convex setting, mixed volumes take the place of intersection numbers.
This can be seen in two steps.

First, $dd^\pm|x|^2=\sum_{i=1}^n dx_i\w d\xi_i$ is the standard Euclidean
super-metric on $\RR^{2n}$. Thus \cite[Example 3.2]{Lagerberg}, \cite[p. 509]{Bern-Arnold}
$$
\frac1{n!}(dd^\pm|x|^2)^n=
dx_1\w d\xi_1\w\cdots\w dx_n\w d\xi_n=
(-1)^{n(n-1)/2}dx_1\w\cdots\w dx_n\w
d\xi_1\w\cdots\w d\xi_n,
$$
that we denote by $(-1)^{n(n-1)/2} dx\w d\xi$ for short.
There are two ways to get rid of the dimensional constant $(-1)^{n(n-1)/2}$.
Berndtsson defines 
\beq
\lb{BerezinModifiedEq}
\int_{\RR^n}d\xi_1\w\cdots\w
d\xi_n=(-1)^{n(n-1)/2}
\eeq
instead of \eqref{BerezinEq}. Lagerberg ignores
the sign altogether which could be considered
either as an oversight or simply as the convention that
forms $dx^\alpha$ commute with forms $d\xi^\beta$, so that
$
dx_1\w d\xi_1\w\cdots\w dx_n\w d\xi_n=
dx_1\w\cdots\w dx_n\w
d\xi_1\w\cdots\w d\xi_n.
$
Either way, 
$$
\bal
|K|=\int_{K} dx=
\int_{K\times\RR^n} dx\w d\xi
=\int_{K\times\RR^n} dx\w d\xi
=\int_{K\times\RR^n} \frac1{n!}(dd^\pm|x|^2)^n,
\eal
$$
and 
assuming $\partial K$ is smooth and 
using Stokes' theorem,
$$
\bal
|K|
&=\int_{\partial K\times\RR^n} d^\pm|x|^2\w\frac1{n!}(dd^\pm|x|^2)^{n-1}
=\frac1n\int_{\partial K\times\RR^n} \sum_{i=1}^nx_i\widehat{dx_i},
\eal
$$
where
$\widehat{dx_i}:=
dx^{\{1,\ldots,i-1,i+1,\ldots,n\}}=
dx^1\w\cdots\w dx^{i-1}\w dx^{i+1}\w\cdots\w dx^n$.
At this point, observe that 
(when $K$ is sufficiently regular)
the gradient of the support function 
of $K$, $\nabla h_K$, 
is a diffeomorphism from $S^{n-1}(1)\subset\RR^n$ to $\partial K$
(Lagerberg claims the same from $\partial K^\circ$ to $\partial K$ \cite[Example 3.2]{Lagerberg}, and that is also correct, observing that $\nabla h_K$
is 0-homogeneous; e.g., 
for $K=[-2,2]$ and $K^\circ=[-1/2,1/2]$, one has $\nabla h_K(x)=2\hbox{\rm sign}(x)$).
Thus, the last integral becomes 
$$
\bal
\frac1n\int_{S^{n-1}(1)\times\RR^n} \sum_{i=1}^n\partial_{x_i}h_K\widehat{d\partial_{x_i}h_K}
&=
\frac1{n!}\int_{S^{n-1}(1)\times\RR^n} d^\pm h_K \wedge (dd^\pm h_K)^{n-1}
\cr
&=
\frac1{n!}\int_{B_2^n\times\RR^n}  (dd^\pm h_K)^{n}
=
\frac1{n!}\int_{\RR^n\times\RR^n}  (dd^\pm h_K)^{n},
\eal
$$
where $B_2^n$ is the unit 2-ball in $\RR^n$ (with $\partial B_2^n=
S^{n-1}(1)$), and the last equality relies on the fact that $(dd^\pm h_k)^{n}$
is actually zero everywhere outside the origin! This can be understood
from the fact that $h_K$ is 1-homogenous, i.e., linear along rays, 
thus its Hessian has a (at least) 1-dimensional kernel everywhere on
$\RR^n\setminus\{0\}$, see \S\ref{toric-HRMA}.




Similar ideas give the following beautiful theorem, a new point of view on the century-old concept of mixed volumes \cite[Proposition 3.12]{Lagerberg}.  

\bthm
Let $K_1,\ldots,K_n\subset \mathcal{K}_0(\RR^n)$ be convex bodies containing the origin in their interior. 
For $f_i\in \mathcal{E}_{K_i}$, the mixed volume equals the super-integral
$$
n!V(K_1,\ldots,K_n)=
\int_{\RR^{2n}}
dd^{\pm}f_1\w
\cdots \w dd^{\pm}f_n.
$$
\ethm





\section{Toric geometry as a bridge}
\lb{toricSec}

Steve could be described as a fanatic fan of toric geometry. He would
test essentially every fact, conjecture, or problem related to
complex geometry that he was interested in at a given moment in the toric world
(e.g., \cite{FZ,RZ1,RZ2,RZ3,SZ1,SZ2,Zel-JSG,Zel,ZZ}). He was also enamoured
with almost everything Donaldson was studying, and Donaldson at that time
was fascinated with toric geometry \cite{Don,Don02,Don-toric}. 
It is then hardly a coincidence that
our first paper \cite{RZ1} was also my first foray 
into the world of toric geometry.
Toric geometry is indeed a very cosmopolitan not-so-little world in its own right, on the 
intersection of algebraic, complex, symplectic, and convex geometry, 
as well as probability, combinatorics, number theory, dynamical systems, and numerical analysis.
It is an excellent area to meet other mathematical nomads and learn
about interesting problems in other fields. This partly explains why Steve
liked toric geometry. 


Concentrating on the convex-complex exchange, objects from complex geometry simplify
drastically on a toric variety to become convex objects associated to Euclidean spaces.
 There is a common misconception that things are `easy' in the presence
of toric symmetry. Indeed, toric varieties are a popular realm to test `general'
conjectures. But perhaps closer to the truth is that there
is a simple trade-off:
the reduction in the number of variables makes it possible
to seek far more granular answers that would be typically out-of-reach
in the 
general (non-toric) complex world.
In other words, oftentimes one can do better than simply resolve 
the special case of a conjecture. Of course, sometimes this can
be misleading as the toric setting might not generalize to the non-toric one.
In sum, research problems concerning toric varieties
can lead to deep mathematical pursuits and insights.
I will try to give some examples for such instances.


\subsection{Toric manifolds: higher-dimensional surfaces of revolution}

The prototype for toric manifold are the simply-connected surfaces of revolution in $\RR^3$.
However, the key feature that generalizes is not the fact that the surface
(topologically  $S^2$) 
embeds  in $\RR^3$ (which is a low-dimensional coincidence) but rather the
underlying (complex!) group action. 
(One can also extend much of what we will discuss here to the genus one case, i.e.,
the torus $S^1\times S^1$ (elliptic curve) that generalizes to abelian 
varieties, another class of manifolds Steve was very fond of, see, e.g., 
\cite{Z-prel},\cite[p. 304]{RZ1},
and the thesis he directed \cite{Feng}; for exposition
and more recent work see \cite{Hult,R-Hult}.

Viewed as a Riemann surface, there is a dense $\CC^*=\CC\setminus\{0\}$ 
inside $\hat\CC:=\CC\cup\{0,\infty\}\cong S^2$, 
and this complex torus is a subgroup of the M\"obius group (of conformal
transformations of $\hat \CC$). If we 
consider geometric objects that are invariant under
the real sub-torus $S^1\subset \CC^*$ then they 
depend solely on a single real variable in $\RR\cong
\CC^*/S^1$. 

Now replace the one-dimensional tori by higher-dimensional ones and you have 
entered into the wonderful world of ``toric geometry." 
This is also a fascinating door into algebraic geometry which is 
encoded by the complex linear group $(\CC^*)^n$. In higher-dimensions
there can be different ways to compactify $(\CC^*)^n$ into a compact complex
manifold. Each such will be a toric variety and the compactification
will be governed by, surprisingly, some elementary but not-at-all-trivial,
geometry of integral lattices and associated convex cones \cite{Fulton}. 


\bdefn
$\bullet\;$ 
a normal variety $X$ \cite[p. 177]{GH} of complex dimension $n$ is called a \emph{toric variety} if it
contains a complex torus $(\CC^*)^n\subset X$ as a dense open subset, together with a
biholomorphic action $(\CC^*)^n \times X\ra X$ of $(\CC^*)^n$ on $X$ that extends 
the natural action of $(\CC^*)^n$ on itself. 

\noindent
$\bullet\;$
A compact K\"ahler manifold $(X,\omega)$ of complex dimension $n$ is a 
\emph{toric \K manifold} if
the compact torus $(S^1)^n$ acts by isometries on $(X,\omega)$ 
and the action extends to a biholomorphic action of a complex torus 
$(\CC^*)^n$ on $X$ with a free, open, dense orbit $X_0\subset X$.
\edefn

A beautiful feature of toric \K manifolds is that they are in one-to-one correspondence
with certain lattice polytopes in $\RR^n$, called {\it Delzant polytopes}.
The polytope can be considered as the moment polytope of the Hamiltonian
$(S^1)^n$-action on $(X,\omega)$ and we refer to \cite{Cannas} for this
point of view. For further exposition see \cite{Apost,JR1,JR2}.

\subsection{The Cauchy problem for the HRMA}
\lb{toric-HRMA}

Not too long after I arrived at Hopkins,
Steve and I started working on the Cauchy problem for the homogeneous
complex \MAeq (HCMA). 
Essentially all previous research on the HCMA involved the Dirichlet
problem: aside from the short-time existence for analytic data that follows 
from the Cauchy-Kowaleskaya theorem, nothing was known about the Cauchy problem.
Thus, it was natural to start with
the simplest setting for the problem, i.e., toric \K manifolds,
where the problem reduces to the homogeneous
real \MAeq (HRMA). Fortunately, almost nothing was known about this problem either
and so we happily set out to~work~on~it. 

Recall from the end of \S\ref{RTSec} that the HCMA corresponds to the equation for 
geodesics in Semmes' geometry on the space $\calH_\omega$ \eqref{HLEq}. The much-studied Dirichlet
problem corresponds to the problem of connecting two metrics by a geodesic, 
very similarly to the interpolation problem. The Cauchy problem corresponds on the other hand to the initial value problem for geodesics emanating from a given metric in a given direction. Its analogue in the setting of interpolation does not seem to have been 
studied. 


The initial value problem for geodesics in the space of toric \K metrics can be written as 
$$
\begin{array}{rrl}
\det\nabla^2\psi
\!\!\!& = & \!\!\!
0, \quad\quad\;\,\mskip2mu \mbox{on} \; [0,T] \times \RR^n,
\cr
\psi(0,\,\cdot\,)
\!\!\!& = &\!\!\!
\psi_0(\,\cdot\,),
\;\; \mbox{on} \;  \RR^n,
\cr
\displaystyle
\frac{\partial\psi}{\partial s}(0,\,\cdot\,)
\!\!\!& = &\!\!\!
\dot\psi_0(\,\cdot\,), \;\; \mbox{on} \; \RR^n,
\end{array} 
$$
with $\nabla^2\psi$ the Hessian in all $n+1$ variables 
$(s,x_1,\ldots,x_n)\in [0,T] \times \RR^n$.
Thus, solutions to HRMA are intuitively expected to have a null-direction
for the Hessian. This is only intuitive, since, as typical, 
one 
extends the meaning of the equation to functions that need not be twice-differentiable.
One of the crucial difficulties is to devise analytic and geometric tools
that make such an intutition rigorous. This can be thought of as one of the conceptual
key points underlying our trilogy \cite{RZ1,RZ2,RZ3}.
Thus, to make the PDE meaningful to study we must
specify ``regularity'' and ``branch of the equation.'' 
Both  of these serve to narrow down the set of possible competitors
and enable to successfully establish an existence and uniqueness theory.
For regularity we assume:

\medskip
\noindent
$\bullet\,$ $\psi_0:\RR^n\ra\RR$ is smooth and strongly convex,

\smallskip
\noindent
$\bullet\, \overline{\nabla_x\psi_0(\RR^n)}=P$, with $P$ the Delzant polytope of $X$.

\smallskip
\noindent
$\bullet\,$ $\dot\psi_0:\RR^n\ra\RR$ is smooth and bounded.

\medskip
The first two assumptions mean that the initial condition corresponds
to a toric \K metric on $X$. Next:

\noindent
\smallskip
\noindent
$\bullet\,$ For the branch we look for a solution $\psi$ that is convex in all $n+1$ variables
and $\overline{\partial_x\psi(s,\,\cdot\,)(\RR^n)}=P$ for each $s\in(0,T]$.

This last assumption means that we require that the curve $s\ra \psi(s,\,\cdot\,)$
corresponds to a path in the space of toric \K metrics on $X$, i.e., that the
\K class is not changing: if it were changing than by the correspondence mentioned in the previous subsection also the moment polytope would be changing. 

For instance, 
$$\psi_0(x)=\log(1+e^{2x})-x, \q \dot\psi_0(x)=\bigg(\frac{e^{2x}-1}{e^{2x}+1}\bigg)^2
$$
satisfies these assumptions with $P=[-1,1]$. The solution to the Cauchy problem
then can be written explicitly \cite[Example 3.1]{RZ2}.

Since this is a PDE problem, the first thing to understand is the weak formulation
or interpretation of the equation. In other words, how can we make sense
of the determinant of the Hessian for a function that is not necessarily twice differentiable? This is an old problem that apparently was solved by Alexandrov
in the first half of the 20th century, under the assumption of convexity.
This is explained in the must-read semi-expository article of Rauch--Taylor \cite{RT}.
Essentially the basic observation is that the fundamental object is not the {\it function} 
$\det\nabla^2 f$ but rather the
{\it measure}
$\det\nabla^2 f\, ds\w dx^1\w\cdots\w dx^{n}
$ 
with that density. The measure is the 
push-forward of the Lebesgue measure on $P$ by the map $\partial\psi^*$.
Here $\psi^*$ denotes the Legendre dual of $\psi$ \cite{rockafellar} and
$\partial$ is the sub-differential operator/map (that is set-valued), generalizing the gradient
to convex but possibly non-differentiable, functions. 
Thus, we may decree that $\psi$ solves the HRMA if the above push-forward
assigns to all open sets in $\RR^n$ zero measure. When $\psi$ is differentiable
but not necessarily twice-differentiable this can be understood geometrically
by duality and means simply that the gradient (in all $n+1$ variables) image of $\psi$
(i.e., the set $\partial\psi([0,T]\times\RR^n)$)
is a measure zero set in $P\times\RR$. There are a number of facts here from convex analysis
that make this all work, but here is the gist. 
The Cauchy problem precisely says that we know the full gradient of $\psi$ 
at time $s=0$: It~is~the~graph 
$$
\{(\nabla_x\psi_0(x),\dot\psi_0(x))\,:\, x\in\RR^n\}.
$$
Duality, via the first variation formula for the Legendre transform, says this is equal
to the graph of $\dot u_0:=-\dot\psi_0\circ(\nabla_x\psi_0)^{-1}$ over $P$!
This is the graph of a nice smooth function over a polytope. Intuitively then, 
the function $\psi$ will solve the HRMA if and only if the full gradient image
remains in that graph (considered as a set in $P\times\RR$). It turns out that for 
some finite time $T_{cvx}$, called the convex lifespan, the image will always satisfy this condition. Thus, the Cauchy problem for the HRMA always has a solution, albeit
for a possibly short time. Moreover, this solution will extend to all time
if and only if $\dot u_0$ is convex on $P$. As soon as it is not convex at a point
there will exist a time such that the gradient image will `fatten up' around that point, 
i.e., accumulate positive mass, and hence $\psi$ will stop solving the equation.
This picture is summarized in \cite[Figure 4, p. 3022]{RZ2}.



As an anecdote, I remember how in April or May 2009 toward the end of the academic year I presented
a draft of \cite{RZ2} in a seminar at Hopkins. The main result I discussed was rather the opposite of what I had just explained! Namely,  in the seminar I presented that 
the IVP 
had a global in time solution for any given $(\psi_0,\dot\psi_0)$. The point of my talk was to show why the 
natural candidate for a solution continued to be a weak solution even when it stopped being smooth (which is exactly the first $s$ for which $\psi_0^* + s\dot u_0$ is no longer
strictly convex on $P$).
Soon thereafter I realized that the opposite was in fact true, it
could not be a solution as soon as it stopped being smooth, or even $C^1$!
I spent a whole academic year (2009--2010) writing up that paper \cite{RZ2} which contained
the result I explained earlier. 
What emerged from that was one of my first real convex-complex exchanges: using
rather intricate tools of convex analysis to prove a result about 
the space of \K metrics. 

Steve and I continued our study of the Cauchy problem for the (general) HCMA. One of the take-aways from the third installment in our series \cite{RZ3} was  
that this sort of phenomenon was special to the toric setting.
That is, without toric symmetry, smooth (or even $C^3$) geodesics in the space of \K metrics may not exist in a given direction {\it even for a short time}
\cite[Theorem 2.7]{RZ3}. We also continued our study of the Cauchy problem for the HRMA
and almost managed to prove a converse to the result from \cite{RZ2}. Namely,
in \cite{RZ2} we showed that the natural candidate solution stops solving the equation
as soon as is stops being smooth. But, perhaps, there is another way to continue it while still solving the equation? In \cite[Theorem 2.14]{RZ3} we showed that there is no $C^1$ way of continuing the solution. The idea was to show that the Cauchy problem for the HRMA is equivalent to a Hamilton--Jacobi equation \cite[Theorem 2.15]{RZ3}.  Since convex functions are merely Lipschitz, it still is a very interesting open problem whether no convex solution exists beyond $T_{cvx}$.

\section{Okounkov bodies and Witt Nystr\"om's Chebyshev transform}

In \S\ref{toricSec} we dwelled a bit about facts that are special
to toric varieties versus facts that generalize from the toric
setting to the general setting. A folklore question for some time
has been, is there a non-toric analogue of a Delzant polytope? 
It turns out the answer is yes, though the polytope is replaced by
a more general convex body, called the Okounkov body of the pair
$(X,L)$ where $X$ is a projective variety and $L$ is an ample
(or even big) line
bundle on it \cite{KavehK,LM}. More recently, this has been extended to non-projective
$X$ as well \cite{trans-OB}. Only countably-many convex bodies
arise as Okounkov bodies on smooth projective manifolds,
and for projective surfaces these are still only polygons, but in higher
dimensions non-polyhedral bodies arise, already for $\PP^2\times\PP^2$
\cite[\S6.3]{LM},\cite[\S3]{KLM}. Already for surfaces, describing the
Okounkov polygon explicitly requires a precise description of the pseudoeffective and nef cones.
For del Pezzo surfaces this has been carried out recently and reveals beautiful
geometric patterns \cite{JR0}.

Given an analogue of the Delzant polytope $P$, one may ask whether more refined
toric objects also have analogues. For instance, the Legendre transform
of an open-orbit \K potential ($u_0:=\psi_0^*$ in the notation of
\S\ref{toric-HRMA}) is a convex function (called the symplectic potential) on $P$ when $X$ is toric.
Witt Nystr\"om showed there is a remarkable way of transforming \K potentials
to convex functions on the Okounkov body for general $X$. He called his 
generalization of the Legendre transform the Chebyshev transform. Several
properties of the symplectic potential carry over to these Chebyshev potentials.
For instance, a beautiful theorem of Witt Nystr\"om shows that the Monge--Amp\`ere
energy can be expressed in terms of the integral of the Chebyshev potential
\cite[Theorem 1.4]{WN},
directly generalizing a  similar fact in the toric world, due to Mabuchi and 
Donaldson \cite{Mab,Don02}. 

In the toric setting the symplectic potentials
themselves are (pointwise) affine along a geodesic in the space of \K metrics.
It was conjectured by R\'eboulet that this carries over to the Chebyshev transform,
but this was disproved recently \cite{JR-Cheby}. The construction therein 
sheds some light on why the linearity holds in the toric setting in the first place:
the toric symmetry forces the flag to be expressed in terms of diagonal matrices
with respect to toric monomoial basis of $H^0(X,L^k)$. Non-toric \K potentials
correspond to more general, non-diagonal, `change-of-basis' matrices. As long as 
these remain uni-lower-diagonal (with respect to the chosen monomial basis) the {\it toric} flag is still preserved
and also the linearity of the Chebyshev 
potentials. However, it is not hard to write down geodesics for which the change-of-basis matrices are not lower-triangular, already on $\PP^1$ \cite[Corollary 1.5]{JR-Cheby}. 

The theory of Okounkov bodies raises many interesting question like these,
concerning analogies to toric geometry. For instance, can the Mabuchi K-energy
and the Ding functional be expressed in terms of the Chebyshev transform?
In a different direction, a discrete analogue of the Okounkov body
has been introduced recently \cite{JR3}, revealing further ways toric geometry and one-dimensional 
phenomena (e.g., Weierstrass gap theory) can be
generalized, raising many interesting questions related to algebraic geometry,
number theory, and~even~microlocal~analysis.


\section{Complex Legendre duality}
\lb{CxLegSeg}

The Legendre transform is one of the main tools of classical convex analysis \cite{rockafellar}. It is a powerful tool that has applications even outside
of mathematics, e.g., in economics, and in mathematics it makes appearances
in areas as diverse  as probability, PDE, symplectic geometry and number theory,
to name a few. Legendre duality is intricately and
intimately related to the real vector space structure of $\RR^n$.
Does it have complex analogues?
It turns out the answer is yes, and moreover, that the classical Legendre
transform itself can be obtained and discovered that way.
The first clue to that comes from a problem raised in the 
celebrated paper of Semmes \cite{S}. Semmes observed that (very roughly) there is 
a certain duality operation on certain real-analytic Lagrangian submanifolds 
in the complexification of the twisted cotangent bundle of a \K manifold, whose
unique self-dual point is a given such 
\cite[p. 543]{S}.

About 30 years later Semmes' intuition was made rigorous in \cite{BCKR2}.
The other starting point is that the only
self-dual function for the Legendre transform in $\RR^n$ is 
$\frac12\sum_{i=1}^nx_i^2$.
The initial idea is that one can decide a priori on a real-analytic function 
and then construct a Legendre-type transform that will have that 
function as its unique self-dual function. For example, assuming
$n=2m$ and starting
with $\phi_{class}(x)=|x|^2/2$ we can consider it in holomorphic variables as 
$\phi_{class}(z)=|z|^2/2$ where $z\in \CC^m$, or equivalently 
$\phi_{class}(z,\bar z)=z\cdot \bar z/2$. 
Then polarizing the variables, and calling $\bar z$ a new variable $w$,
we may consider $\phi_{class}(z,w)$ as a function on $\CC^n=\CC^m\times \CC^m$.
The clasical Legendre transform is then simply (up to harmless factors of $2$
and identification of $\CC^m$ with $\RR^n$)
$
f^*(w)=\sup_z[2\hbox{\rm Re}\phi_{class}(z,w)-f(z)].
$
Thus, if we replace $|x|^2/2$ with any real-analytic $\phi$ we formally obtain 
a new transform. There are a number of technical issues with this (e.g., 
over what domain is the supremum being taken? Will $\phi$ be 
its unique fixed point?), but that is the gist.  
In fact, this same construction translates to any {\it compact} \K manifold,
yielding such transforms that have absolutely nothing to do with the vector
space structure that seemed so inherently central to the classical Legendre
transform. For the \K setting a crucial necessary technical tool is Calabi's
diastasis function that allows to prove that $\phi$ is the unique fixed point.
These complex Legendre transforms turn out to be (local) isometries of the space of \K metrics \cite{BCKR2}, in fact the only ones \cite{Lempert-local},
just as in Semmes' vision, and closely related to Lempert's 
beautiful construction of symmetries of the HCMA \cite{Lempert} (see also \cite{Darvas-JEMS}).


\section{A convex-complex approach to the Mahler conjectures}\label{sec2}

Let $K\subset \R^n$ be a {convex body}, i.e., a compact convex set with non-empty interior and recall the definition of its polar $K^\circ$ from \S\ref{BBLSec}.  
We call $K$  symmetric
when  $K=-K\defeq \{x\in\R^n: -x\in K\}.$
Denote by 
$    B_p^n\defeq \{x\in\R^n: {x_1^p+\ldots+ x_n^p}\leq 1\}, 
$ the standard closed $p$-ball, with $B_\infty^n:=[-1,1]^n$.

Let $A\in GL(n,\R)$. 
While $AK$ and $K$ could have wildly differing volume
(with respect to the Lesbegue measure $d\lambda$), the volume of the product body
$AK\times (AK)^\circ\subset \R^n\times\R^n$ is equal
to that of $K\times K^\circ$.
This leads to the following $GL(n,\R)$-invariant functional on convex bodies \cite[p. 95]{mahler}.
\begin{definition}
\label{Mahlerdef}
The \textit{Mahler volume} of $K$ is
$    \M(K)\defeq n! |K\times K^\circ|= n! |K||K^\circ|.
$\end{definition}

Crude bounds on $\M$ were demonstrated by Mahler already in 1939 \cite[(6)]{mahler}.
In 1987,
Bourgain--Milman showed that there exists an unspecified but uniform $c>0$ independent of $K$ and $n$ such that 
\cite[Corollary 6.1]{bourgain-milman} 
\begin{equation}\label{bm_int}
    \M(K)\geq c^n,
\end{equation}
that---aside from determining the best value of $c$---is optimal in an
asymptotic sense \cite[pp. 149--150]{ryabogin-zvavitch}. 
Mahler's first conjecture asserts that $c$ should be $4$ if $K$ is symmetric, attained by both $B_1^n$ and $B_\infty^n$ \cite[p. 96]{mahler}.
Mahler's second conjecture asserts that
among general convex $K$,
\begin{equation}\label{bm_int_ns}
    \M(K)\geq \frac{(n+1)^{n+1}}{n!}, 
\end{equation}
attained by simplices centered at the origin
(this is, by Stirling's formula, asymptotic to $e^n$)
\cite[(1)]{mahler2} \cite[p. 564]{schneider}.

The best known constant to date is $c=\pi$ for the first conjecture
(and $c=4$ for $n=3$ \cite{iriyeh-shibata})
and $\pi/2$ (and $c=2$ for $n=3$) for the second; a beautiful complex proof using contour
integrals is due to Berndtsson \cite{berndtsson1} (cf. \cite{kuperberg2}).

More than a decade ago, Nazarov proposed a complex approach to the first Mahler
conjecture using the Bergman kernel of 
the `tube domains 
over $K$',
$T_K\defeq 
\R^n+ \sqrt{-1}(\mathrm{int}\,K)
\subset \C^n
$ \cite{nazarov}.
By applying H\"ormander's technique for $L^2$ estimates on solutions to
the $\bar\partial$-equation \cite{hormander2} he obtained the first complex proof of
\eqref{bm_int} for symmetric $K$ with $c=\pi^3/16$ (this was
extended to general $K$ with $c=\pi/4$ \cite{MR}). Shortly afterwards
Blocki gave an alternative proof of one of Nazarov's estimates \cite{blocki2,blocki}.
There remained hope that perhaps this approach of Nazarov could
give a line of attack to Mahler's first conjecture, but recently
it was realized that Nazarov's approach cannot yield $c=4$. 

The key observation is that the Bergman kernel in Nazarov's approach
can be realized as an {\it $L^1$ version of the Mahler volume} \cite{MR,BMR}.
Thus, the H\"ormander--Nazarov approach can, at best, yield a proof
of an $L^1$ version of the Mahler conjecture, that is also equivalent
to a conjecture of Blocki on Bergman kernels~of~tube~domains.

However, the new observation we wish to point out here is a new
{\it convex-complex} approach to both Mahler's conjectures. The approach
relies on reducing the Mahler conjecture (either version) to {\it two}
separate conjectures, one that can be attacked using purely complex
methods, and another that can be attacked using purely convex methods.

The $L^p$-Mahler volume $\M_p$ is derived by realizing that the classical
support function (that underlies Legendre duality!) has $L^p$ counterparts
\cite{MR,BMR}.

\bdefn
    For $p\in(0,\infty)$, 
$            h_{p,K}(y)\defeq \log\left[\int_K e^{p\langle x,y\rangle}\frac{dx}{|K|}\right]^{\frac1p}, \; y\in\R^n
$, and
\begin{equation*}
            \M_p(K)\defeq |K|\int_{\R^n}e^{-h_{p,K}(y)}dy.
        \end{equation*}
    
\edefn


Berndtsson--Mastrantonis--Rubinstein posed the $L^p$-Mahler conjecture (and, similarly,
a version for general $K$) \cite{BMR}.

\begin{conjecture}
\lb{MpConj}
For symmetric $K\subset \RR^n$, 
$
\M_p(K)\ge        \M_p(B_\infty^n).
$
    
\end{conjecture}

We pose: 
    \begin{conjecture}
    \lb{MpMConj}
        For symmetric $K\subset\R^n$, 
        \begin{equation*}
            \frac{\M}{\M_1}\big(K\big)\geq 
            \frac{\M}{\M_1}\big(B_\infty^n\big) = \left( \frac{4}{\pi^2}\right)^n. 
        \end{equation*}
    \end{conjecture}
    A similar conjecture can be made (up to translation) in the non-symmetric case with the simplex as the minimizer. 
If this last conjecture holds, it would bridge the gap between the inequality 
\begin{equation}
\lb{NazarovEq}
    \M(K)\geq \frac{\M_1(K)}{4^n}
\end{equation}
that is implicit in Nazarov's work and Conjecture \ref{MpConj} for $p=1$,
namely, $\M_1(K)\ge \pi^{2n}$. Indeed, combining Conjecture \ref{MpConj}
and \eqref{NazarovEq} only yields \eqref{bm_int} with $c=\pi^2/4$ while
combining both Conjectures \ref{MpConj} and \ref{MpMConj} would yield
Mahler's conjectured $c=4$. Of course, this could all be considered
as an approach without substance, except that in fact
it does work in the case $n=2$, where Conjecture 
\ref{MpConj} was obtained in \cite{MR2} and Conjecture \ref{MpMConj}
is expected \cite{MR3}. It could, of course, be a futile approach as Conjecture
\ref{MpMConj} could very well turn-out to be as difficult to prove
as Mahler's conjecture itself, even if a complex proof is found for 
Conjecture~\ref{MpConj}~with~$p=1$.

\section{Epilogue}

I first met Steve Zelditch in Summer 2007 during a summer-long visit to the Technion
while still a graduate student. 
This meeting sparked a rather intense collaboration over the next five years
with our first paper appearing already in early 2008.
I obtained my Ph.D. from M.I.T. in 2008, and was offered tenure-track positions
at Johns Hopkins and Stony Brook.  
Having just started to learn microlocal analysis through many discussions with Steve, I thought it would be a good idea to spend some more time to continue learning it from one of its leaders,
and I decided to do a post-doctoral period before settling more permanently.
So I turned down the offers and chose instead
to spend the year 2008--2009 at Johns Hopkins (which was also Steve's last
year in Johns Hopkins before moving to Northwestern) before continuing to Stanford 
for three more years. I learned a lot from Steve during those years as well as in endless discussions in years to come. When Steve left us in 2022 I realized just
how different collaborating with Steve was from any other collaboration I had.
Steve Zelditch was a bigger than life proponent of interactions between convex and complex 
(among many other things). 
He pioneered several such. His mathematical choices strongly influenced my own,
and continue to inspire many. In this article I shared some of the mathematical
journey that I have been through in the lands of convex and complex, and where
I continue to wa(o)nder. 
For more personal recollections on Steve
I refer the reader to \cite{R23-notices}. 

Steve, thank you for being my teacher.

\bibliographystyle{amsplain}

\end{document}